\documentclass[11pt]{amsart}
\usepackage{psfig}
\usepackage{amssymb}
\usepackage{amscd}

\newcommand{\sph}[1]{{\mathbb S}^{#1}}
\newcommand{\hype}[1]{{\mathbb H}^{#1}}
\newcommand{\euc}[1]{{\mathbb E}^{#1}}
\newcommand{\ints}{{\mathbb Z}}

\newcommand{\reals}{{\mathbb R}}
\newcommand{\complex}{{\mathbb C}}

\newcommand{\slice}{M}

\newcommand{\mink}[1]{{{\mathbb R}_1^{#1}}}

\theoremstyle{plain}
\newtheorem{thm}{Theorem}[section]

\newtheorem{lemma}[thm]{Lemma}
\newtheorem{cor}[thm]{Corollary}

\theoremstyle{definition}

\theoremstyle{remark}

\newtheorem*{case0}{Case 0}
\newtheorem*{case1}{Case 1}
\newtheorem*{case2}{Case 2}
\newtheorem*{case3}{Case 3}

\begin{document}
\title[Spacelike slices]{$3$-manifolds which are spacelike slices of flat spacetimes}
\author{Kevin P. Scannell} 
\address{Department of Mathematics and Computer Science, Saint Louis University, St. Louis, MO 63103} 
\email{scannell@slu.edu}
\date{October 16, 2000}
\thanks{The author was partially supported by NSF grant DMS-0072515.}

\begin{abstract}
\noindent
We continue work initiated in a 1990 preprint of Mess
giving a geometric parameterization of the moduli space 
of classical solutions to Einstein's equations in $(2+1)$ dimensions 
with cosmological constant $\Lambda = 0$ or $-1$ (the case $\Lambda = +1$
has been worked out in the interim by the present author).
In this paper we make a first step toward the $(3+1)$-dimensional case
by determining exactly which closed $3$-manifolds $M^3$ arise as 
spacelike slices of flat spacetimes, and by finding all
possible holonomy homomorphisms $\pi_1(M^3) \to ISO(3,1)$.
\end{abstract}

\maketitle

\section{Introduction}

This paper has as its starting point the work of Witten on
$(2+1)$-dimensional gravity from the late 1980's \cite{Wi89},
\cite{Wi89a}, \cite{Wi91}.  In the first of these papers
he asks for a description of the ``space of all classical
solutions'' to Einstein's equations in $(2+1)$-dimensions.  
In a remarkable paper shortly thereafter \cite{Me90}, Geoffrey Mess 
interprets and completely
solves this problem for flat and anti-de~Sitter 
spacetimes, giving a
geometric parameterization roughly in terms of the ``degeneration''
of spacelike slices.   The de~Sitter case
was handled in \cite{Sc96} and \cite{Sc99}, where a new phenomenon appears:
an infinite family of solutions with identical holonomy representations
(Witten remarks on the possible physical significance of this in 
\cite[\S 6]{Wi91}).
This paper represents a first step in extending these
results to the case of constant curvature $(3+1)$-dimensional
spacetimes.   

Throughout this paper, $\slice$ will denote a closed, 
connected $3$-manifold.  Our current understanding of $3$-manifold
topology owes a great deal to Thurston's introduction
of geometric techniques in the 1970's.
In particular, his Geometrization Conjecture \cite{Sc84}, \cite{Th82}
says that given a closed $3$-manifold, there is a canonical process 
by which it can be cut open so that the 
resulting pieces are {\em geometric}:  this means they can be given
Riemannian metrics which are locally isometric to one of eight 
simply connected Riemannian homogeneous spaces.   
As some of these ``model spaces'' will arise in the statement 
of the main theorem
and in the course of its proof, it is worth setting up some notation for them.
There are, of course, the three constant curvature model spaces: 
Euclidean space $\euc{3}$, hyperbolic space $\hype{3}$, 
and the $3$-sphere $\sph{3}$;
two product spaces $\hype{2} \times \reals$ and $\sph{2} \times \reals$;
and three $3$-dimensional Lie groups $Nil$, $Solv$, and 
$\widetilde{SL_2\reals}$ equipped with natural left-invariant metrics.
If $X$ is one of these eight spaces, we say a $3$-manifold is 
{\em modelled on $X$} if it admits a metric locally isometric to $X$.
When a manifold $M$ is modelled on $\hype{3}$, we often say simply
that $M$ is {\em hyperbolic}.   If we exclude $\hype{3}$ and $Solv$,
then $M$ is modelled on one of the six remaining spaces if and only
if $M$ is a {\em Seifert fiber space} \cite[Thm. 5.3]{Sc84}.

We say $M$ is a {\em spacelike slice of a flat spacetime} if 
there is a $(3+1)$-dimensional Lorentzian manifold $N$
locally isometric to Minkowski space 
$\mink{4}$ and an embedding $f : M \hookrightarrow N$
such that $f(M)$ is spacelike and has trivial normal bundle.
One immediate consequence of our main theorem is that a spacelike slice 
of a flat spacetime is geometric, and in fact can only be modelled on
three of the eight model spaces:
\begin{thm} \label{T:flatmain}
$\slice$ is a spacelike slice of a flat spacetime
if and only if $\slice$ is modelled on $\hype{3}$, $\euc{3}$, or 
$\hype{2} \times \reals$.
\end{thm}
The manifolds modelled on $\euc{3}$ and $\hype{2} \times \reals$ 
are understood completely.  There are exactly ten closed $3$-manifolds
modelled on $\euc{3}$ \cite[\S 3.5]{Wo84} all of which are finitely covered by
the $3$-torus $T^3 = \sph{1} \times \sph{1} \times \sph{1}$.   
Though there are infinitely many  
$3$-manifolds modelled on $\hype{2} \times \reals$, each is finitely
covered by $\Sigma \times \sph{1}$ where $\Sigma$ is a closed 
surface of negative Euler characteristic (see Lemma \ref{L:covers} below).    
Thus we have as a corollary:
\begin{cor} \label{C:flatmain}
If $M$ is a spacelike slice of a flat spacetime and is not hyperbolic,
then a finite cover of $M$ is homeomorphic to $\Sigma \times \sph{1}$,
where $\Sigma$ is a closed, orientable surface, $\Sigma \ncong \sph{2}$.
\end{cor}
The class of hyperbolic $3$-manifolds is much richer than
the others by the work of Thurston which shows that ``most'' $3$-manifolds
are hyperbolic.    In fact, the Geometrization
Conjecture predicts that any closed, irreducible $3$-manifold
with infinite fundamental group not containing a $\ints \oplus \ints$
subgroup is hyperbolic.

Mess observes \cite[p. 55]{Me90} that the $(2+1)$-dimensional
case generalizes in part to the $(3+1)$-dimensional case; in particular
he claims the following result which will follow immediately from our
proof of Theorem \ref{T:flatmain}:
\begin{cor} \label{C:messclaim}
Let $M$ be a spacelike slice of a flat spacetime.
If the linear holonomy representation $L : \pi_1(M) \to O(3,1)$ 
is irreducible, then $M$ is hyperbolic.
\end{cor}
Also noteworthy is the work of H. Waelbroeck \cite{Wa95} which 
is closely related to this paper and \cite{Me90}.  
Flat spacetimes homeomorphic to $M \times \reals$ 
are studied in terms of solutions of the so-called $B \wedge F$ theory.
Our main result shows however
that many of the solutions found in \cite{Wa95}
(e.g. for the $Nil$, $Solv$, and $\widetilde{SL_2\reals}$ cases)
do not correspond to actual flat spacetimes with spacelike slices.

For comparison with Theorem \ref{T:flatmain}, 
we recall the corresponding statement from the de~Sitter case
(particular examples illustrating this case have appeared in the
physics literature; see \cite{BH99}, \cite{Ez94}, \cite{MW93}).
A manifold is {\em conformally flat} if it admits a 
(locally) conformally flat Riemannian metric; for $3$-manifolds
this is equivalent \cite{Sc96} to the existence of a {\em flat conformal} or 
{\em M\"obius} structure.
\begin{thm} \label{T:deSmain}
\cite{Sc96}
$\slice$ is a spacelike slice of a de~Sitter spacetime
if and only if $\slice$ is conformally flat.
\end{thm}
Together these results show that there are many more possibilities for
slices of de~Sitter spacetimes than for flat spacetimes:
for example, all manifolds modelled on the constant curvature or product 
model spaces are conformally flat.  
Also notable is the fact that the connected sum of conformally flat manifolds 
is conformally flat \cite{Ku78b}
and even some manifolds modelled on 
$\widetilde{SL_2\reals}$ are
conformally flat \cite{GLT88},
providing large classes 
of examples not present in the flat case.
Unfortunately, no
classification of conformally flat $3$-manifolds is known in general,
even assuming the Geometrization Conjecture.   The best results in
this direction are due to M. Kapovich; see for instance \cite{Ka93}.

It should be emphasized that the statement of the main theorem
is purely topological --
it says nothing about how a given $\slice$ arises as a spacelike
slice, its induced Riemannian metric, or about the moduli space 
of spacetimes having a given $\slice$ as a Cauchy
surface.   In the de~Sitter case, these questions were worked out completely
in \cite{Sc99}: the moduli space of de~Sitter domains of dependence
$\slice \times \reals$ is identified with an appropriate deformation space
of conformally flat metrics on $\slice$, thus reducing the question 
to a widely studied problem in Riemannian geometry.
The classification theorems in \cite{Me90} and \cite{Sc99} basically come
from a convexity result for the causal horizon of a 
spacelike slice, followed by a careful analysis of the 
geometric structure of the causal horizon.
Our goal in writing this paper, in contrast, was to derive 
as much as possible purely from results in $3$-manifold topology,
hopefully returning to a study of the causal horizon in a sequel.

An added bit of information which falls out in the course of the
proof of Theorem \ref{T:flatmain}  is the determination
of all possible holonomy homomorphisms $\pi_1(\slice) \to ISO(3,1)$.  
Of course this falls well short of describing the moduli
space of flat spacetimes, as it is not {\it a priori} true
that constant curvature spacetimes $\slice \times \reals$
are parameterized by their holonomy homomorphisms
(indeed this is false in the de~Sitter case as mentioned above).
Furthermore, one needs to beware of certain holonomy homomorphisms 
which are inadmissible because
they only arise from spacetimes homeomorphic to $\slice \times \reals$
where the slices $\slice \times \{ t \}$ are not spacelike.

The author wishes to thank Lars Andersson, 
Bill Goldman, Geoff Mess and Steve Harris for
useful conversations related to this work.

\section{Basic results} \label{S:basic}
Let $X$ be a Riemannian or Lorentzian homogeneous space,
and let $G$ be its isometry group.
If $V$ is a smooth manifold of the same dimension as $X$,
then we can define a {\em $(G,X)$-structure} on $V$ to be
a maximal atlas of coordinate charts 
$\{\phi_\alpha : U_\alpha \to X\}$ on $V$
such that the transition functions are 
given by the action of elements of $G$.
From this data, a standard argument
constructs  a $(G,X)$-structure on the universal cover $\widetilde{V}$,
a {\em developing map} $\mathcal{D} : \widetilde{V} \to X$,
and a {\em holonomy homomorphism} $\phi : \pi_1(V) \to G$
satisfying the following equivariance condition
$$
\mathcal{D}(\gamma \cdot x) = \phi(\gamma) \cdot \mathcal{D}(x)
$$
for all $\gamma \in \pi_1(V)$ and $x \in \widetilde{V}$
(see \cite{Go88a} for a nice discussion of these notions).
One can show that the existence of a $(G,X)$-structure
is equivalent to the existence of a Riemannian or Lorentzian
metric everywhere locally isometric to the model space $X$;
in particular, one could rephrase our discussion of the Geometrization
Conjecture in this language.

Let $\slice \hookrightarrow N$ be a spacelike slice of a flat spacetime. 
Since $M$ has trivial normal bundle, we might as well assume
$N$ is homeomorphic to $M \times (0,1)$.
The ideas just introduced provide us with a developing map 
$\mathcal{D} : \widetilde{N} \to \mink{4}$
and a holonomy homomorphism $\phi : \pi_1(N) \to ISO(3,1)$
satisfying the equivariance condition above.
Here $ISO(3,1)$ denotes the full isometry group of $\mink{4}$
which we will describe in detail momentarily.
We compose these functions with the inclusions
$\widetilde{\slice} \hookrightarrow \widetilde{N}$ 
and $i_\ast : \pi_1(\slice) \stackrel{\cong}{\to} \pi_1(N)$
respectively to obtain
a spacelike immersion 
$\mathcal{D} : \widetilde{\slice} \to \mink{4}$
and a {\em holonomy group} $\Gamma = \phi(i_\ast \pi_1(M)) \subset ISO(3,1)$.

Our first result is that this immersion is actually an embedding
and that $\widetilde{\slice} \cong \reals^3$.
The proof given here is special to Minkowski
space (the analogous result is false for de~Sitter space).    
Similar theorems are obtained for general classes of spacetimes
by Harris in \cite{Ha87}, \cite{Ha88}.
\begin{lemma} \label{L:embed}
The image of $\mathcal{D} : \widetilde{\slice} \to \mink{4}$ is a graph over 
$\euc{3} = \{ \mathbf{v} \in \mink{4} \mid v_4 = 0 \}$.  Any two points of
$\mathcal{D}(\widetilde{\slice})$ are spacelike separated.
\end{lemma}

\begin{proof}
The obvious projection $\mink{4} \to \euc{3}$ restricts to a 
distance-increasing local diffeomorphism 
$P : \mathcal{D}(\widetilde{\slice}) \to \euc{3}$ 
since $\mathcal{D}(\widetilde{\slice})$ is everywhere spacelike.  
As $M$ is closed, the induced Riemannian metric is complete, 
and thus $P$ must be proper.  This makes $P$ a covering map, hence
a global diffeomorphism.  Finally $\mathcal{D}(\widetilde{\slice})$ 
is a graph because $\widetilde{M}$ is connected.

If $\mathbf{p}$ and $\mathbf{q}$ 
are two points of the image which are null or timelike
separated, consider the path between them given by the intersection of 
$\mathcal{D}(\widetilde{\slice})$ and an indefinite two-dimensional
plane containing $\mathbf{p}$ and $\mathbf{q}$.  
The ``secant line'' joining $\mathbf{p}$ and $\mathbf{q}$ 
in this plane has slope greater than $1$, so the mean value theorem 
implies that some tangent vector to this path is null or timelike,
a contradiction.
\end{proof}
In all that follows we will identify $\widetilde{M}$ with its image
in $\mink{4}$.

Every isometry of $\mink{4}$
can be written uniquely as $\mathbf{x} \mapsto A\mathbf{x} + \mathbf{b}$,
where the {\em linear part} $A$ lies in $O(3,1)$
and $\mathbf{b} \in \mink{4}$. 
If we let $L: ISO(3,1) \to O(3,1)$ be the homomorphism projecting
to the linear part, we have the following short exact sequence
$$
1 \to \mink{4} \to ISO(3,1) \stackrel{L}{\to} O(3,1) \to 1.
$$
Given a subgroup $\Gamma$ of $ISO(3,1)$,
define $T(\Gamma) = \ker L \bigr|_\Gamma$;
we call $T(\Gamma)$ the {\em translational subgroup} of $\Gamma$.
There is a corresponding short exact sequence for $\Gamma$
$$
1 \to T(\Gamma) \to \Gamma \stackrel{L}{\to} L(\Gamma) \to 1
$$
which is central to our study of possible holonomy groups.

\begin{lemma} \label{L:disc}
Let $\Gamma \subset ISO(3,1)$ be the holonomy group of a spacelike slice.
Then:
\begin{enumerate}
\item $\Gamma$ is a discrete, torsion-free subgroup of $ISO(3,1)$.
\item $T(\Gamma)$ consists of spacelike vectors and is isomorphic to $\ints^k$,
for $k = 0,1,2,$ or $3$.
\item $L(\Gamma)$ leaves invariant the spacelike subspace spanned 
by $T(\Gamma)$.
\end{enumerate}
\end{lemma}

\begin{proof}
The first part is straightforward by Lemma \ref{L:embed} -- 
since $\widetilde{M} \cong \reals^3$, $\Gamma$ has finite cohomological
dimension and is therefore torsion-free.  Discreteness follows because
$\Gamma$ acts properly discontinuously on $\widetilde{M}$.
Because all pairs of points in $\widetilde{M}$ are spacelike-separated, 
it is clear that $T(\Gamma)$ consists only of spacelike vectors.  It is 
isomorphic to $\ints^k$ since it is a discrete subgroup of $\mink{4}$
and $k$ cannot be $4$ or else $T(\Gamma)$ would have to contain a 
non-spacelike vector.
For the last part, if $\gamma_1 \in \Gamma$ is
given by $\mathbf{x} \mapsto A\mathbf{x} + \mathbf{b}$
and $\gamma_2 \in T(\Gamma)$ is translation by $\mathbf{t}$,
then $\gamma_1 \gamma_2 \gamma_1^{-1}$ is easily computed to be
the element $\mathbf{x} \mapsto \mathbf{x} + A \mathbf{t}$
of $T(\Gamma)$.
\end{proof}

\begin{lemma} \label{L:Aus}
If $\Gamma$ is a discrete subgroup of $ISO(3,1)$ with $L(\Gamma)$ 
indiscrete, then $L(\Gamma)$ is virtually solvable.
\end{lemma}

\begin{proof}
A theorem of Auslander (see \cite[Thm. 8.24]{Ra72})
says that $\overline{L(\Gamma)}^0$ is always solvable
and is non-trivial since $L(\Gamma)$ is indiscrete. 
The closed, connected, solvable, non-trivial Lie subgroups of $SO(3,1)^0$
are easy to write down as in \cite{CG74}; in particular, the
set of points $F$ on the sphere at infinity $\partial \hype{3}$ fixed by 
$\overline{L(\Gamma)}^0$ must consist of one or two points.
The stabilizer in $SO(3,1)^0$ of a point at infinity
is isomorphic to the group $Sim^+(\reals^2)$ of orientation-preserving
similarities of $\reals^2$, which is a solvable group.
Since $L(\Gamma)$ normalizes $\overline{L(\Gamma)}^0$,
it leaves $F$ invariant and therefore has
a subgroup of index at most two fixing $F$ pointwise,
and therefore conjugate into $Sim^+(\reals^2)$.
The lemma follows.
\end{proof}

The remaining lemmas are well-known results from 
$3$-manifold topology.   We record the ones which will be
used repeatedly in \S \ref{S:main}.

\begin{lemma} (see \cite{He76}) \label{L:Wald}
Suppose $\widetilde{\slice} \cong \reals^3$ and 
$\ints^3 \subseteq \pi_1(\slice)$.
Then $\slice$ is finitely covered by the $3$-torus $T^3$.
\end{lemma}

\begin{proof}
While a much more general result is proved 
in \cite{He76}, we give a simpler proof sufficient for our needs.
Consider the cover
$\widehat{\slice}$ of $\slice$ corresponding to the $\ints^3$
subgroup.   Since $\widetilde{\slice} \cong \reals^3$, 
$\widehat{\slice}$ is a $K(\ints^3, 1)$ hence is homotopy equivalent
to $T^3$.  This implies that 
$H_3(\widehat{\slice}) \cong H_3(T^3) \cong \ints$, thus
$\widehat{\slice}$ is closed and the covering 
$\widehat{\slice} \to \slice$ is finite.
Finally, appealing to \cite{Wa68}, we have 
$\widehat{\slice} \cong T^3$.
\end{proof}

The next result explains the topological structure
of manifolds modelled on $\euc{3}$ and $\hype{2} \times \reals$.
It connects Theorem \ref{T:flatmain} and Corollary \ref{C:flatmain}
and will also be used in \S \ref{S:main}.

\begin{lemma} (see \cite{Sc84}) \label{L:covers}
$M$ is modelled on $\euc{3}$ if and only if it is finitely covered by
the $3$-torus, and $M$ is modelled on $\hype{2} \times \reals$ 
if and only if it is finitely covered by $\Sigma \times \sph{1}$, 
where $\Sigma$ is a closed, orientable surface of genus at least two.
\end{lemma}


\section{Proof of main theorem} \label{S:main}
Before embarking on the proof of Theorem \ref{T:flatmain}
we note that its main content is
the ``only if'' part, which excludes many kinds
of $3$-manifolds from being spacelike slices of flat spacetimes.
In particular, manifolds modelled on five of
Thurston's eight geometries cannot be spacelike slices;
it is useful in traversing the proof to keep some of
these exclusions in mind.  For instance, it follows 
immediately from Lemma \ref{L:embed} that manifolds
modelled on $\sph{3}$ or $\sph{2} \times \reals$ 
do not arise.  An important element
of the proof of the main theorem is to show that
the Euler number of a Seifert fiber space which 
is a spacelike slice must be zero, excluding manifolds
modelled on $Nil$ or $\widetilde{SL_2\reals}$.
$Solv$ and $Nil$ are interesting since it is 
possible to find flat spacetimes homeomorphic to 
$\slice \times \reals$ where $\slice$ is modelled on $Solv$ or $Nil$,
but the main theorem says that the slices $\slice \times \{ t \}$ 
can never be spacelike. 

\begin{proof} (of Theorem \ref{T:flatmain})
The ``if'' half is easy as each of the spaces $\hype{3}$,
$\euc{3}$, and $\hype{2} \times \reals$ embeds in $\mink{4}$
(though there can be geometrically distinct embeddings as 
we will see for $\euc{3}$).
We have 
\begin{itemize}
\item $\hype{3} =  \{ \mathbf{v} \in \mink{4} \mid \langle \mathbf{v}, \mathbf{v} \rangle = -1, v_4 > 0 \}$
\item $\euc{3} =  \{ \mathbf{v} \in \mink{4} \mid v_4 = 0 \}$
\item $\hype{2} \times \reals = \{ \mathbf{v} \in \mink{4} \mid v_2^2 + v_3^2 - v_4^2 = -1, v_4 > 0 \}$.
\end{itemize}
If $\slice$ is modelled on a model space $X$, then it
can be realized as $\slice \cong X / \Gamma$, where $\Gamma$ is 
a discrete, cocompact subgroup of the isometry group $Isom(X)$.   
It is easy to see that if $X$ is one of the three examples above embedded 
in $\mink{4}$, then there is a corresponding embedding of its isometry 
group in $ISO(3,1)$, and that
any discrete subgroup of $Isom(X)$ acts discontinuously
on a regular neighborhood of $X$ in $\mink{4}$.
The quotient of a small regular neighborhood is therefore a flat spacetime
containing $M \cong X / \Gamma$ as a spacelike slice, as desired.

For the ``only if'' half, let 
$\slice$ be a spacelike slice of a flat spacetime with 
holonomy group $\Gamma \subset ISO(3,1)$.
The proof is broken down into four cases, depending on the rank
of the translational subgroup $T(\Gamma)$  (Lemma \ref{L:disc}).
Note that the cases get easier as we go along, because the
presence of a large normal abelian subgroup of $\pi_1(M)$ 
typically puts strong topological constraints on $M$.

\begin{case0}
Suppose $T(\Gamma) = 0$.  This means that $L$ injects $\Gamma$
into $O(3,1)$, i.e. $\Gamma \cong L(\Gamma)$.
Our assumption that spacelike slices have
trivial normal bundles means that $L(\Gamma)$ actually lies in
the orthochronous subgroup $O_\uparrow(3,1)$ which coincides
with the full isometry group of $\hype{3}$.  
Now $L(\Gamma)$ is either discrete or indiscrete.
If it is discrete, then it is also 
cocompact for cohomological reasons.  Since $M$ is aspherical, 
it is homotopy equivalent to the closed 
hyperbolic $3$-manifold $\hype{3} / L(\Gamma)$ 
and a result of Gabai-Meyerhoff-Thurston \cite{GMT96} implies that $M$
is itself hyperbolic.  In fact, it will turn out that this is
the only possibility for the holonomy when $M$ is hyperbolic.
Thus in all remaining cases we will be proving that $M$ is 
modelled on $\euc{3}$ or $\hype{2} \times \reals$.
In light of Lemma \ref{L:covers}, it suffices to show that 
$M$ has a finite cover homeomorphic to $\Sigma \times \sph{1}$,
where $\Sigma$ is a closed, orientable surface, $\Sigma \ncong \sph{2}$.
We will exploit this fact in all that follows by freely passing
to finite covers of $\slice$ without changing notation.
Also note that $L(\Gamma)$ is reducible in all remaining cases,
yielding Corollary \ref{C:messclaim}.

If $L(\Gamma)$ is indiscrete,
Lemma \ref{L:Aus} shows that $L(\Gamma) \cong \Gamma$ is 
virtually solvable.
By the main result of \cite{EM72}, we may pass to a finite cover
and assume $M$
is a torus bundle over $\sph{1}$ with monodromy 
$\theta : \pi_1(T^2) \to \pi_1(T^2)$
represented by a matrix
in $SL(2,\ints)$, by abuse of notation also denoted $\theta$.
Clearly $\theta$ has finite order if and only if $M$ is finitely covered
by a $3$-torus, we will assume $\theta$ has infinite order.
Let $t \in \Gamma$ denote an element inducing
the monodromy, i.e. $t x t^{-1} = \theta(x)$ for all $x \in \pi_1(T^2)$.
The image of the fiber subgroup 
$\pi_1(T^2) \cong \ints \oplus \ints$ of $\Gamma$ under $L$ 
must consist either of elements leaving invariant 
a geodesic in $\hype{3}$ (generated by loxodromics or irrational
elliptics) or of parabolics with a common fixed point at infinity.
In the first case,
since $t$ normalizes $\pi_1(T^2)$,
$L(t)$ must leave the geodesic invariant
(indeed it must fix it pointwise since $L(t)$ has infinite order). 
But this implies that $t$ commutes with $\pi_1(T^2)$ 
which means $\theta$ is the identity and $\Gamma \cong L(\Gamma) \cong \ints^3$.
This contradicts the hypothesis that $\theta$ has infinite order,
or alternatively shows directly that $M \cong T^3$ in this case by
the proof of Lemma \ref{L:Wald}.

In the second case, identify 
$\partial \hype{3} \cong \complex \cup \{ \infty \}$
and conjugate so that the parabolics fix $\infty$.
If $\theta$ has a $1$-eigenvalue,
let $g$ denote an element of
$\pi_1(T^2) \cong \ints \oplus \ints$ such that $\theta(g) = g$
and write $L(g)$ as $z \mapsto z + z_0$ for some $z_0 \neq 0 \in \complex$.
Since $t$ normalizes $\pi_1(T^2)$,
$L(t)$ must also fix $\infty$; write it as $z \mapsto az + b$
for some $a,b \in \complex$.   
But then with this notation, the relation $[L(t),L(g)] = 1$
becomes $a z_0 = z_0$.   
Thus $a = 1$, $L(t)$ is parabolic,
$t$ commutes with $\pi_1(T^2)$, and we are done as above. 
The final possibility is that $\theta$ is of infinite order and has
no $1$-eigenvalue,
in which case it has two
real eigenvalues of absolute value not equal to one.
In this case, there is, up to conjugacy, only one
possibility for $L(\Gamma)$, and we use this in the Appendix
to show that $\Gamma$ cannot be discrete, contradicting
Lemma \ref{L:disc} (1).

This completes Case $0$.
\end{case0}

\begin{case1}
Suppose $T(\Gamma) \cong \ints$.   We may assume $M$ is orientable
by passing to a double cover.  Though we can get away with
less, we might as well use the Seifert Fiber Space Theorem
of Mess \cite{Me}, Gabai \cite{Ga92}, and Casson-Jungreis \cite{CJ94}
which states that a closed, orientable, irreducible $3$-manifold
whose fundamental group contains a normal $\ints$ is a Seifert fiber
space.  The proof of this result amounts to showing that the
quotient group $\pi_1(M) / \ints$ (namely $L(\Gamma)$ in our setup) 
is the fundamental group of a closed $2$-orbifold $B$, 
which by passing to a finite cover, we
can assume is a closed, orientable surface of genus $g \geq 1$.
It follows that $M$ is modelled on $\euc{3}$, $\hype{2} \times \reals$, $Nil$,
or $\widetilde{SL_2\reals}$.   Excluding the final two possibilities
amounts to showing that the Euler number $e$ of $M$ is zero.  

Our original proof that $e = 0$ amounted to finding the place
where the Euler number appears in the Hochschild-Serre sequence 
for $H^1(\Gamma, \mink{4})$;
we have chosen to give a more geometric argument here.
The group $\Gamma$ has a presentation
of the form \cite[p. 91]{Or72}
$$
\Gamma = \langle a_1, b_1, \ldots, a_g, b_g, h \mid [a_i, h] = [b_i, h] = 1,
h^e = [a_1,b_1] \cdots [a_g,b_g] \rangle;
$$
the notation is meant to be obvious: $h$ generates the normal 
subgroup $T(\Gamma) \cong \ints$
and the other generators project to $L(\Gamma)$, 
the fundamental group of the base $B$.
Write $\mathbf{x} \mapsto \mathbf{x} + \mathbf{v}$
for the translation $h$, and
$\mathbf{x} \mapsto L(a_i) \mathbf{x} + t(a_i)$,
$\mathbf{x} \mapsto L(b_i) \mathbf{x} + t(b_i)$
for the other generators.
If we let $W = (\mink{4})^{L(\Gamma)}$ be the subspace left
invariant by $L(\Gamma)$,
then the commutation relations imply that
$L(a_i) \mathbf{v} = \mathbf{v}$ and $L(b_i) \mathbf{v} = \mathbf{v}$,
and so $\mathbf{v} \in W$.
Now focus on the final relation: the left-hand side is the translation 
$\mathbf{x} \mapsto \mathbf{x} + e \mathbf{v}$
while the right-hand side's translational part
is a combination of the $t(a_i)$ and $t(b_i)$:
$$
e \mathbf{v} = \sum_{i = 1}^{g} L(c_{i-1})(I-L(a_i b_i a_i^{-1})) t(a_i) + 
L(c_{i-1} a_i)(I - L(b_i a_i^{-1} b_i^{-1})) t(b_i)
$$
where $c_j = [a_1, b_1] \cdots [a_j, b_j]$.
The linear map being applied to $t(a_i)$ (resp. $t(b_i)$) 
in the expression above is called the {\em Fox derivative}
$\frac{\partial R}{\partial a_i}$ (resp.
$\frac{\partial R}{\partial b_i}$)
of the usual surface group relator $R$ (see \cite{Fo53}, \cite{Go84}). 
These partial derivatives can be combined neatly into a single 
{\em Fox differential} $dR : (\mink{4})^{2g} \to \mink{4}$, 
in which case
$$
e \mathbf{v} = dR(t(a_1), t(b_1), \ldots, t(a_g), t(b_g)).
$$
In \cite[\S 3.7]{Go84}, Goldman 
has a nice argument showing 
that the image of $dR$ is exactly $W^\perp$.
Thus $\mathbf{v} \in W$ and $e \mathbf{v} \in W^\perp$, so
$$
0 = \langle \mathbf{v}, e \mathbf{v} \rangle = 
e \langle \mathbf{v}, \mathbf{v} \rangle.
$$
Since $\mathbf{v}$ is spacelike (this is essential -- see \cite{AM59}),
$\langle \mathbf{v}, \mathbf{v} \rangle \neq 0$, and so $e = 0$.

This completes Case $1$.
\end{case1}

\begin{case2}
Suppose $T(\Gamma) \cong \ints^2$.  
A theorem in Hempel \cite[Thm. 11.1]{He76} shows that
$L(\Gamma)$ has two ends and therefore
has a finite index subgroup isomorphic to $\ints$. 
As usual, we will pass to a finite cover without changing notation
and assume $L(\Gamma) \cong \ints$.  Stallings' theorem \cite{St62}
implies that $M$ fibers overs the circle with torus fibers.
Let $A \in O(3,1)$ be a generator of $L(\Gamma)$.
It leaves invariant the spacelike $\euc{2}$ spanned by 
$T(\Gamma) \subset \mink{4}$, and in fact, acts by linear isometries.
Since $A$ normalizes the lattice $T(\Gamma)$, the action of $A$ 
on $\euc{2}$ must have finite order.   It follows that $M$ is finitely
covered by the $3$-torus, and we conclude that $M$ is
modelled on $\euc{3}$.

This completes Case $2$.
\end{case2}

\begin{case3}
Suppose $T(\Gamma) \cong \ints^3$.  
Lemma \ref{L:Wald} implies that $M$ is finitely covered by the $3$-torus
and is therefore modelled on $\euc{3}$.
\end{case3}
This completes the proof of the main theorem.
\end{proof}

\appendix
\section{}
This appendix contains the proof of a result used
in Case 0 of the proof of the main theorem.
It is relegated to an appendix because it relies on a
technical cohomology calculation which would have 
interrupted the flow of the exposition to an unacceptable degree.

The goal is to dispose of a particular class of solvable groups $\Gamma$
which are fundamental groups of torus bundles over $\sph{1}$ with
``hyperbolic'' monodromy.  We will retain all of the notations
used in proof of the main theorem.
In particular, the monodromy of the fibration 
$\theta : \pi_1(T^2) \to \pi_1(T^2)$
is represented by a matrix
$A \in SL(2,\ints)$
which has two real eigenvalues, say $\lambda > 1$ and $1/\lambda < 1$.
Let $x_1$ and $x_2$ be the standard generators of $\pi_1(T^2)$
and write the corresponding parabolics $L(x_j)$ as $z \mapsto z + w_j$, 
$j = 1,2$.
As in \S 3, we write $L(t)$ as
$z \mapsto az + b$ for some $a,b \in \complex$.   
The relation $txt^{-1} = \theta(x)$ applied to the two 
generators collates into the following matrix equation:
$$
A
\left(
\begin{matrix}
w_1\\
w_2  
\end{matrix}
\right)
=
\left(
\begin{matrix}
a w_1\\
a w_2  
\end{matrix}
\right).
$$
In other words, $a$ must be an eigenvalue of $\theta$ (say, $a = \lambda$),
and the $w_j$ are the components of the corresponding eigenvector
(in particular, they are real and irrationally related).
For simplicity, we choose $w_1 = 1$.

\renewcommand{\thethm}{A.1}
\begin{lemma}
With $L(\Gamma)$ given as above, $\Gamma$ must be indiscrete.
\end{lemma}

\begin{proof}
In fact, we will show that the subgroup $\pi_1(T^2) \subset \Gamma$
must be indiscrete.   The idea of the proof is to view $\Gamma$ 
as obtained from
the indiscrete group $L(\Gamma)$ by adding translations, and to show
that there is no way of doing so which makes $\pi_1(T^2)$ discrete.
Such a choice of translations is a cocycle in $H^1(\Gamma, \mink{4})$;
that is, a function $c : \Gamma \to \mink{4}$ satisfying
the {\em cocycle relation}
$$
c(gh) = c(g) + L(g) \cdot c(h)
$$
for all $g,h \in \Gamma$.
The coboundaries (change of basepoint) are cocycles of the form
$$
c(g) = (1 - L(g)) v
$$
for some fixed $v \in \mink{4}$.
There is, of course, a restriction map to the fiber subgroup
$$
H^1(\Gamma, \mink{4}) \to H^1(\pi_1(T^2), \mink{4})^{\langle t \rangle},
$$
where the notation is meant to indicate that restricted classes
are invariant under the action of $t$:
\begin{equation} \label{E:invariance}
[t \cdot c] = [c] \text{ where } (t \cdot c)(x) := L(t) c(\theta^{-1}(x)).
\end{equation}
The lemma will follow from the fact that the group
$H^1(\pi_1(T^2), \mink{4})^{\langle t \rangle}$ vanishes.
The cocycle relation applied to $x_1 x_2 = x_2 x_1$ 
means that for any cocycle 
$c \in H^1(\pi_1(T^2), \mink{4})$ we must have
\begin{equation} \label{E:cocycle}
(1 - L(x_2)) c(x_1) = (1 - L(x_1)) c(x_2).
\end{equation}
Equation (\ref{E:invariance}) (really two equations for $x_1$ and $x_2$)
and Equation (\ref{E:cocycle}) are linear in 
$c(x_1)$ and $c(x_2)$ and can be solved explicitly 
by choosing a basis for $\mink{4}$ and elements of
$O(3,1)$ representing the generators of $L(\Gamma)$.  For instance, following
\cite{CG74}, we can choose the first two basis elements to be null
vectors fixed by $L(t)$ (the second fixed by $L(x_j)$)
and the last two basis elements to be spacelike (the final one
also fixed by $L(x_j)$).
This yields:
$$
L(x_1) = 
\left(
\begin{matrix}
1 & 0 & 0 & 0\\
\frac{1}{2} & 1 & 1 & 0\\
1 & 0 & 1 & 0\\
0 & 0 & 0 & 1
\end{matrix}
\right),
\hspace{1ex}
L(x_2) = 
\left(
\begin{matrix}
1 & 0 & 0 & 0\\
\frac{w_2^2}{2} & 1 & w_2 & 0\\
w_2 & 0 & 1 & 0\\
0 & 0 & 0 & 1
\end{matrix}
\right),
\hspace{1ex}
L(t) = 
\left(
\begin{matrix}
\lambda & 0 & 0 & 0\\
0 & \lambda^{-1} & 0 & 0\\
0 & 0 & 1 & 0\\
0 & 0 & 0 & 1
\end{matrix}
\right).
$$
Some linear algebra shows that there is a one-dimensional
solution set to Equations
(\ref{E:invariance}) and (\ref{E:cocycle});
namely by taking $c(x_1) = (0,1,0,0)$ and $c(x_2) = (0,w_2,0,0)$.
This solution is a coboundary however, since 
$c(x_j) = (1 - L(x_j)) v$ for 
$v = (0,0,-1,0)$ as a simple calculation shows.
Thus $H^1(\pi_1(T^2), \mink{4})^{\langle t \rangle} = 0$,
proving the lemma.
\end{proof}

\bibliographystyle{amsplain}
\bibliography{data}

\providecommand{\bysame}{\leavevmode\hbox to3em{\hrulefill}\thinspace}
\begin{thebibliography}{10}

\bibitem{AM59}
L.~Auslander and L.~Markus, \emph{Flat {L}orentz $3$-manifolds}, Mem. Amer.
  Math. Soc. \textbf{30} (1959), 1--60.

\bibitem{BH99}
I.~Bengtsson and S.~Holst, \emph{de~{S}itter space and spatial topology},
  Classical Quantum Gravity \textbf{16} (1999), no.~11, 3735--3748.

\bibitem{CJ94}
A.~J. Casson and D.~S. Jungreis, \emph{Convergence groups and {S}eifert fibered
  $3$-manifolds}, Invent. Math. \textbf{118} (1994), no.~3, 441--456.

\bibitem{CG74}
S.~S. Chen and L.~Greenberg, \emph{Hyperbolic spaces}, Contributions to
  analysis (L.~V. Ahlfors, I.~Kra, B.~Maskit, and L.~Nirenberg, eds.), Academic
  Press, New York-London, 1974, pp.~49--87.

\bibitem{EM72}
B.~D. Evans and L.~E. Moser, \emph{Solvable fundamental groups of compact
  $3$-manifolds}, Trans. Amer. Math. Soc. \textbf{168} (1972), 189--210.

\bibitem{Ez94}
K.~Ezawa, \emph{Classical and quantum evolutions of the de {S}itter and the
  anti-de {S}itter universes in $2+1$ dimensions}, Phys. Rev. D \textbf{49}
  (1994), no.~10, 5211--5226.

\bibitem{Fo53}
R.~H. Fox, \emph{Free differential calculus. {I}. {D}erivation in the free
  group ring}, Ann. of Math. (2) \textbf{57} (1953), no.~3, 547--560.

\bibitem{Ga92}
D.~Gabai, \emph{Convergence groups are {F}uchsian groups}, Ann. of Math. (2)
  \textbf{136} (1992), no.~3, 447--510.

\bibitem{GMT96}
D.~Gabai, R.~Meyerhoff, and N.~Thurston, \emph{Homotopy hyperbolic
  $3$-manifolds are hyperbolic}, MSRI Preprint 1996-058, 1996.

\bibitem{Go84}
W.~M. Goldman, \emph{The symplectic nature of fundamental groups of surfaces},
  Adv. Math. \textbf{54} (1984), no.~2, 200--225.

\bibitem{Go88a}
\bysame, \emph{Geometric structures on manifolds and varieties of
  representations}, Geometry of {G}roup {R}epresentations (W.~M. Goldman and
  A.~R. Magid, eds.), Contemp. Math., vol.~74, Amer. Math. Soc., Providence,
  1988, pp.~169--197.

\bibitem{GLT88}
M.~L. Gromov, H.~B.~Lawson Jr., and W.~P. Thurston, \emph{Hyperbolic
  $4$-manifolds and conformally flat $3$-manifolds}, Inst. Hautes \'Etudes Sci.
  Publ. Math. \textbf{68} (1988), 27--45.

\bibitem{Ha87}
S.~G. Harris, \emph{Complete codimension-one spacelike immersions}, Classical
  Quantum Gravity \textbf{4} (1987), no.~6, 1577--1585.

\bibitem{Ha88}
\bysame, \emph{Complete spacelike immersions with topology}, Classical Quantum
  Gravity \textbf{5} (1988), no.~6, 833--838.

\bibitem{He76}
J.~P. Hempel, \emph{$3$-manifolds}, Ann. of Math. Stud., vol.~86, Princeton
  Univ. Press, Princeton, 1976.

\bibitem{Ka93}
M.~\`E. Kapovich, \emph{Flat conformal structures on $3$-manifolds, {I}:
  {U}niformization of closed {S}eifert manifolds}, J. Differential Geom.
  \textbf{38} (1993), 191--215.

\bibitem{Ku78b}
R.~S. Kulkarni, \emph{On the principle of uniformization}, J. Differential
  Geom. \textbf{13} (1978), no.~1, 109--138.

\bibitem{Me90}
G.~Mess, \emph{Lorentz spacetimes of constant curvature}, MSRI Preprint
  90-05808, 1990.

\bibitem{Me}
\bysame, \emph{The {S}eifert conjecture and groups which are coarse
  quasiisometric to planes}, Preprint, 1990.

\bibitem{MW93}
J.~Morrow-Jones and D.~M. Witt, \emph{Inflationary initial data for generic
  spatial topology}, Phys. Rev. D \textbf{48} (1993), no.~6, 2516--2528.

\bibitem{Or72}
P.~Orlik, \emph{Seifert manifolds}, Lecture Notes in Math., vol. 291,
  Springer-Verlag, New York-Berlin-Heidelberg, 1972.

\bibitem{Ra72}
M.~S. Raghunathan, \emph{Discrete subgroups of {L}ie groups}, Ergeb. Math.
  Grenzgeb., vol.~68, Springer-Verlag, New York-Berlin-Heidelberg, 1972.

\bibitem{Sc96}
K.~P. Scannell, \emph{Flat conformal structures and causality in de {S}itter
  manifolds}, Ph.D. thesis, University of California, Los Angeles, 1996.

\bibitem{Sc99}
\bysame, \emph{Flat conformal structures and the classification of de {S}itter
  manifolds}, Comm. Anal. Geom. \textbf{7} (1999), no.~2, 325--345.

\bibitem{Sc84}
G.~P. Scott, \emph{The geometries of $3$-manifolds}, Bull. London Math. Soc.
  \textbf{15} (1984), no.~5, 401--487.

\bibitem{St62}
J.~R. Stallings, \emph{On fibering certain $3$-manifolds}, Topology of
  3-manifolds and related topics, Prentice-Hall, Inc., Englewood Cliffs, 1962,
  pp.~95--100.

\bibitem{Th82}
W.~P. Thurston, \emph{$3$-dimensional manifolds, {K}leinian groups and
  hyperbolic geometry}, Bull. Amer. Math. Soc. (N.S.) \textbf{6} (1982), no.~3,
  357--381.

\bibitem{Wa95}
H.~Waelbroeck, \emph{${B}\wedge {F}$ theory and flat spacetimes}, Comm. Math.
  Phys. \textbf{170} (1995), no.~1, 63--78.

\bibitem{Wa68}
F.~Waldhausen, \emph{On irreducible $3$-manifolds which are sufficiently
  large}, Ann. of Math. (2) \textbf{87} (1968), no.~2, 56--88.

\bibitem{Wi89}
E.~Witten, \emph{$2+1$ dimensional gravity as an exactly soluble system},
  Nuclear Phys. B \textbf{311} (1989), no.~1, 46--78.

\bibitem{Wi89a}
\bysame, \emph{Topology-changing amplitudes in $2+1$ dimensional gravity},
  Nuclear Phys. B \textbf{323} (1989), no.~1, 113--140.

\bibitem{Wi91}
\bysame, \emph{Quantization of {C}hern-{S}imons gauge theory with complex gauge
  group}, Comm. Math. Phys. \textbf{137} (1991), no.~1, 29--66.

\bibitem{Wo84}
J.~A. Wolf, \emph{Spaces of constant curvature}, 5 ed., Publish or Perish,
  Inc., Houston, 1984.

\end{thebibliography}
\end{document}